\documentclass[11pt,twoside]{amsart}
\usepackage{amsmath, amsthm, amscd, amsfonts, amssymb, graphicx, color, alltt}
\usepackage[bookmarksnumbered, plainpages]{hyperref}

 \usepackage[usenames,dvipsnames]{pstricks}
 \usepackage{epsfig}
 \usepackage{pst-grad} % For gradients
 \usepackage{pst-plot} % For axes

\usepackage[utf8]{inputenc}
\usepackage{pstricks-add}

\usepackage{enumerate}% package command enumerate
\usepackage{subfig} %  two table in side together

\usepackage{ragged2e}
\usepackage{alltt}
\newtheorem{thm}{Theorem}[section]

\newtheorem{lem}[thm]{Lemma}
\newtheorem{prop}[thm]{Proposition}

\newtheorem{exa}[thm]{Example}
\newtheorem{nota}[thm]{Notation}

\numberwithin{equation}{section}

\newtheorem{remark}{Remark}
\newcommand{\gyr}[2]{{\mathrm{gyr}[{#1}]}{#2}}

\newcommand{\Aut}[1]{\mathrm{Aut}\,{(#1)}}

\title{\bf Construction of Gyrogroups of Order $2^n$ by Cyclic 2-Groups}

\author[A. A. Ungar]{Abraham A. Ungar}
\address{Department of Mathematics, North Dakota State University, Fargo, ND 58108-6050, USA}
\email{abraham.ungar@ndsu.edu}
\thanks{}
\author[M. A. Salahshour]{Mohammad Ali Salahshour}
\address{Department of Mathematics, Savadkooh Branch, Islamic Azad University, Savadkooh, Iran}
\email{salahshour@iausk.ac.ir}
\thanks{}
\author[K. Mavaddat Nezhaad]{Kurosh Mavaddat Nezhaad$^\star$}
\address{\justifying{Department of Pure Mathematics, Faculty of Mathematical Sciences, University of Kashan, Kashan
		87317$-$53153, Iran}}
\email{kuroshmavaddat@gmail.com}
\thanks{$^\star$Corresponding author (Email:kuroshmavaddat@gmail.com)}
\date{}
\begin{document}

\maketitle
\begin{center}
\dedicatory{This work is in memory of Professor Ali Reza Ashrafi, a forbearing wise teacher.}
\end{center}
\begin{abstract}
A gyrogroup is a structure constituting from a non-empty set and a binary operation such that satisfying the left identity, and left inverse conditions, and also has the associative-like law said to be left gyroassociativity and left loop property. In this paper, we propose a method  for constructing new gyrogroups of order $2^n$ by a cyclic 2-group that is $\mathbb{Z}_{2^n}$, where $n\geq 3$.

\vskip 3mm


\vskip 3mm


\end{abstract}

\bigskip

\section{Introduction}
\noindent A pair $(G,\oplus)$ consisting of a non-empty set $G$ and a binary operation $\oplus$ on $G$ is called a groupoid.
Let $(G,\oplus)$ be a groupoid. A bijection from $G$ to itself is called an automorphism of $G$ if $\phi(a\oplus b)=\phi(a)\oplus\phi(b)$ for all $a,b\in G$. The set of all automorphisms of $G$ is denoted by $\Aut{G,\oplus}$. It is easy to see that $\Aut{G,\oplus}$ forms a group under function composition.
A groupoid $(G,\oplus)$ is called a \textit{gyrogroup} if the following conditions hold:
\begin{enumerate}

	\item [(G1)] there exists an element $0 \in G$ such that for all $x \in G$, $0 \oplus x = x$;
	
	\item [(G2)] for each $a \in  G$, there exists $b \in G$ such that $b \oplus a = 0$;
	
	\item [(G3)] there exists a function $\mbox{gyr}: G \times G \longrightarrow \Aut{G,\oplus}$ such that for every 
	$a, b, c \in G$, $a \oplus (b \oplus c) =  (a \oplus b) \oplus \gyr{a,b}{c}$, where  $\gyr{a,b}{c} = \mathrm{gyr}(a,b)(c)$;
	
	\item [(G4)] for each $a, b \in G$, $\gyr{a,b}{} = \gyr{a \oplus b, b}{}$.

\end{enumerate}

Note that these axioms imply their right counterparts. It is easy to see that a group is a gyrogroup if we define gyroautomorphisms to be the identity automorphism. For every $a,b \in G$, the mapping $\gyr{a,b}{}$ is called the gyroautomorphism generated by $a$ and $b$. 
The gyrogroup structure is a result of a pioneering work of Abraham Ungar in the study of Lorentz group \cite{6,7}. A gyrogroup $G$ is called a gyrocommutative if  and only if for all $a,b\in G$,
$$a\oplus b=\gyr{a,b}{(b\oplus a)}$$
Throughout this paper, our notations are standard and can be taken mainly from \cite{1,8}. We refer interested readers to consult the survey \cite{5} for a complete account of the history of gyrogroup. We also refer to \cite{42,44} for subgyrogroups, gyrogroup homomorphisms, and quotient gyrogroup. Our calculations are done with the aid of GAP \cite{45}. For more Information on classification of gyrogroups of small orders of orders at most 31 up to Isomorphism except for 24, 27, and 30 readers can see the \cite{11}. In this paper, we make a gyrogroup of order $2^n$ by the cyclic 2-group $\mathbb{Z}_{2^n}$ for $n\geq 3$ that is denoted by $G_2(n)$. Ashrafi et al. in \cite{12} constructed a gyrogroup $G(n)$ of order $2^n$ such that in this paper we denote it by $G_1(n)$.
The elements of $\mathbb{Z}_{2^n}$ are as follows:
$$\mathbb{Z}_{2^n}=\{0,1,2,\cdots,2^{n-1}-1,2^{n-1},2^{n-1}+1,\cdots,2^n-1\}.$$

%%%%%%%%%%%%%%%%%%%%%%%%%%%%%%%%%
\begin{nota}\label{not1}
{\em
Let $P(n)=\{0,1,2,\cdots,2^{n-1}-1\}$, $H(n)=\{2^{n-1},2^{n-1}+1,\cdots,2^n-1\}$ and $G_2(n)=P(n)\cup H(n)$. It is clear $P(n)\cong \mathbb{Z}_{2^{n-1}}$ and $H(n)=P(n)+m$ that $m=2^{n-1}$. Hence $G_2(n)=P(n)\cup (P(n)+m)$.
Also, let
$O_P = \{i\in P(n)\ \mid \ \text{$i$ is odd}\}$,
$E_P = \{i\in P(n)\ \mid \ \text{$i$ is even}\}$,
$O_H = \{i\in H(n)\ \mid \ \text{$i$ is odd}\}$ and
$E_H = \{i\in H(n)\ \mid \ \text{$i$ is even}\}$.
Define the binary operation $\oplus$ on $G_2(n)$ as follows:
$$i\oplus j =
\begin{cases}
t & :(i,j)\in (P(n)\times P(n))\cup\Bigl[(H(n)\times H(n))-(E_H\times O_H)\Bigr]\\[2mm]
t+m & :(i,j)\in (P(n)\times H(n))\cup\Bigl[(H(n)\times P(n))-(E_H\times O_P)\Bigr]\\[1mm]
s & :(i,j)\in E_H\times O_H \\[1mm]
s+m & :(i,j)\in E_H\times O_P
\end{cases}$$
which $t,s\in P(n)$ and
\[\left\{\begin{split}
t & \equiv i+j \pmod m\\
s & \equiv i+j+\frac{m}{2} \pmod m
\end{split}
\right.\]
Obviously, the operation $\oplus$ is well defined.
}
\end{nota}
Where the greatest common divisor of positive integers $r$ and $s$ is denoted by $(r,s)$.
We start this section by the following simple lemma:
%%%%%%%%%%%%%%%%%%%%%%%%%%%%%%%%%
\begin{lem}[\cite{12}]\label{lem2}
Assume notation \ref{not1}.
\begin{enumerate}[$I)$]
\item $(\frac{m}{2}-1,m)=(\frac{m}{2}+1,m)=1$.
\item  Suppose $x\equiv y\pmod m$. If $x,y\in P(n)$ or $x,y\in H(n)$, then $x=y$.
\end{enumerate}
\end{lem}
%%%%%%%%%%%%%%%%%%%%%%%%%%%%%%%%%
%%%%%%%%%%%%%%%%%%%%%%%%%%%%%%%%%%
\begin{lem} \label{lem1}
Assume notation \ref{not1} and $A(i)$ be a map on $G_2(n)$ as following,
\[A(i) =
\begin{cases}
r & :i\in O_P \\
r+m & :i\in O_H \\
i & :\text{Otherwise}
\end{cases}\]
where $r\in P(n)$ and $r\equiv i+\displaystyle\frac{m}{2} \pmod m$. Then $A(i)\in \Aut{G_2(n),\oplus}$.
\end{lem}

\begin{proof}
Clearly $A(i)$ is a well defined and a bijection map on $G_2(n)$. It is enough to show that $A(i)$ is a homomorphism on $(G_2(n),\oplus)$. To do this, assume $i,j\in G_2(n)$ be two arbitrary elements. We consider four separate cases as follows:
\begin{enumerate}
%%%%%%%%%%%%%%%%%%%%%%%%%%%%%%%%%%%%%%%%%%%%%55
\item $i,j\in P(n)$. It is the first case we consider in four steps as follows:
\begin{enumerate}
%*****************
\item $i,j \in O_P$. By definition $\oplus$, $i\oplus j=t_1\in E_P$ such that $t_1\equiv i+j\pmod m$.
By assumption
\begin{equation}\label{eq1}
A(i\oplus j)=A(t_1)=t_1\equiv i+j\pmod m
\end{equation}
Also, $A(i)=r_1$, $A(j)=r_2$ such that
\[r_1\equiv i+\frac{m}{2}\pmod m \hspace*{1cm} \text{and} \hspace*{1cm} r_2\equiv j+\frac{m}{2}\pmod m. \]
Since $r_1,r_2\in O_P$ then by the last equations and the definition $\oplus$,
\begin{equation}\label{eq2}
A(i)\oplus A(j)=r_1\oplus r_2\equiv r_1+r_2\equiv i+j\pmod m
\end{equation}
Since $A(i)\oplus A(j),A(i\oplus j)\in P(n)$ by equations \ref{eq1}, \ref{eq2} and lemma \ref{lem2}
$A(i\oplus j)=A(i)\oplus A(j)$.
%*****************
\item $i,j \in E_P$. By definition $\oplus$, $i\oplus j=t_1\in E_P$ such that $t_1\equiv i+j\pmod m$.
By assumption
\begin{equation}\label{eq3}
A(i\oplus j)=A(t_1)=t_1\equiv i+j\pmod m
\end{equation}
Also, $A(i)=i$, $A(j)=j$ hence by definition $\oplus$,
\begin{equation}\label{eq4}
A(i)\oplus A(j)=i\oplus j\equiv i+j\pmod m
\end{equation}
Since $A(i)\oplus A(j),A(i\oplus j)\in P(n)$ by equations \ref{eq3}, \ref{eq4} and lemma \ref{lem2}
$A(i\oplus j)=A(i)\oplus A(j)$.
%*****************
\item $i\in O_P$ and $j\in E_P$. By definition $\oplus$,
$i\oplus j=t_1\in O_P$ such that $t_1\equiv i+j\pmod m$. By assumption
\begin{equation}\label{eq5}
A(i\oplus j)=A(t_1)=r_1\equiv t_1+\frac{m}{2}\equiv i+j+\frac{m}{2}\pmod m
\end{equation}
Also $A(j)=j$ and $A(i)=r_1$  such that $r_1\equiv i+\frac{m}{2}\pmod m$. By definition $\oplus$,
\begin{equation}\label{eq6}
A(i)\oplus A(j)=r_1\oplus j\equiv r_1+j\equiv i+\frac{m}{2}+j\pmod m
\end{equation}
Since $A(i)\oplus A(j),A(i\oplus j)\in P(n)$ by equations \ref{eq5}, \ref{eq6} and lemma \ref{lem2}
$A(i\oplus j)=A(i)\oplus A(j)$.
\item $i\in E_P$ and $j\in O_P$. The proof is similar to (c) and so it is omitted.
\end{enumerate}
%%%%%%%%%%%%%%%%%%%%%%%%%%%%%%%%%%%%%
\item $i,j\in H(n)$. It is the second case we consider in four steps as follows:
\begin{enumerate}
%*****************
\item $i,j \in O_H$. By definition $\oplus$, $i\oplus j=t_1\in E_P$ such that $t_1\equiv i+j\pmod m$.
By assumption
\begin{equation}\label{eq7}
A(i\oplus j)=A(t_1)=t_1\equiv i+j\pmod m
\end{equation}
Also, $A(i)=r_1+m$, $A(j)=r_2+m$ such that:
\[r_1\equiv i+\frac{m}{2}\pmod m \hspace*{1cm} \text{and} \hspace*{1cm} r_2\equiv j+\frac{m}{2}\pmod m. \]
Since $r_1,r_2\in O_P$ then by the last equations and the definition $\oplus$,
\begin{equation}\label{eq8}
A(i)\oplus A(j)=(r_1+m)\oplus (r_2+m)\equiv r_1+r_2\equiv i+j\pmod m
\end{equation}
Since $A(i)\oplus A(j),A(i\oplus j)\in P(n)$ by equations \ref{eq7}, \ref{eq8} and lemma \ref{lem2}
$A(i\oplus j)=A(i)\oplus A(j)$.
%*****************
\item $i,j \in E_H$. The proof is similar to (1)(b)  and so it is omitted.
%By definition $\oplus$, $i\oplus j=t_1\in E_H$ such that $t_1\equiv i+j\pmod m$.
%By assumption
%\begin{equation}\label{eq3}
%A(i\oplus j)=A(t_1)=t_1\equiv i+j\pmod m
%\end{equation}
%Also, $A(i)=i$, $A(j)=j$ hence by definition $\oplus$,
%\begin{equation}\label{eq4}
%A(i)\oplus A(j)=i\oplus j\equiv i+j\pmod m
%\end{equation}
%Since $A(i)\oplus A(j),A(i\oplus j)\in P$ by Equations \ref{eq3} and \ref{eq4} $A(i\oplus j)=A(i)\oplus A(j)$.
%*****************
\item $i\in O_H$ and $j\in E_H$. By definition $\oplus$, $i\oplus j=t_1\in O_P$ such that $t_1\equiv i+j\pmod m$.
By assumption
\begin{equation}\label{eq9}
A(i\oplus j)=A(t_1)=r_1\equiv t_1+\frac{m}{2}\equiv i+j+\frac{m}{2}\pmod m
\end{equation}
Also $A(j)=j$ and $A(i)=r_1+m\in O_H$  such that $r_1\equiv i+\frac{m}{2}\pmod m$. By definition $\oplus$,
\begin{equation}\label{eq10}
A(i)\oplus A(j)=(r_1+m)\oplus j\equiv (r_1+m)+j\equiv i+\frac{m}{2}+j\pmod m
\end{equation}
Since $A(i)\oplus A(j),A(i\oplus j)\in P(n)$ by equations \ref{eq9}, \ref{eq10} and lemma \ref{lem2} $A(i\oplus j)=A(i)\oplus A(j)$.
%*****************
\item $i\in E_H$ and $j\in O_H$. By definition $\oplus$, $i\oplus j=s\in O_P$ such that
 $s\equiv i+j+\frac{m}{2}\pmod m$.
 By assumption
\begin{equation}\label{eq11}
A(i\oplus j)=A(s)=r\equiv s+\frac{m}{2}\equiv i+j\pmod m
\end{equation}
Also $A(i)=i$ and $A(j)=r_1+m\in O_H$ such that $r_1\equiv j+\frac{m}{2}\pmod m$. By definition $\oplus$,
\begin{equation}\label{eq12}
A(i)\oplus A(j)= i\oplus(r_1+m) \equiv i+(r_1+m)+\frac{m}{2}\equiv i+j\pmod m
\end{equation}
Since $A(i)\oplus A(j),A(i\oplus j)\in P(n)$ by equations \ref{eq11}, \ref{eq12} and lemma \ref{lem2} $A(i\oplus j)=A(i)\oplus A(j)$.
\end{enumerate}
%%%%%%%%%%%%%%%%%%%%%%%%%%%%%%%%%%%%%
\item $i\in P(n)$ and $j\in H(n)$. It is the third case we consider four in steps as follows:
\begin{enumerate}
%*****************
\item $i \in O_P$ and $j \in O_H$. By definition $\oplus$, $i\oplus j=t+m\in E_H$ such that $t\in P(n)$ and $t\equiv i+j\pmod m$.
By assumption
\begin{equation}\label{eq13}
A(i\oplus j)=A(t+m)=t+m\equiv  i+j\pmod m
\end{equation}
Also $A(i)=r_1\in O_P$ and $A(j)=r_2+m\in O_H$  such that
\[r_1\equiv i+\frac{m}{2}\pmod m \hspace*{1cm} \text{and} \hspace*{1cm} r_2\equiv j+\frac{m}{2}\pmod m. \]
By definition $\oplus$,
\begin{equation}\label{eq14}
A(i)\oplus A(j)=r_1\oplus (r_2+m)\equiv r_1+(r_2+m)\equiv i+j\pmod m
\end{equation}
Since $A(i)\oplus A(j),A(i\oplus j)\in H(n)$ by equations \ref{eq13}, \ref{eq14} and lemma \ref{lem2} $A(i\oplus j)=A(i)\oplus A(j)$.
%*****************
\item $i \in E_P$ and $j \in E_H$. It is clear $A(i)=i$ and $A(j)=j$. By definition $\oplus$, $i\oplus j=t+m\in E_H$ such that $t\in P(n)$ and $t\equiv i+j\pmod m$. Therefore
\[A(i\oplus j)=A(t+m)=t+m=i\oplus j=A(i)\oplus A(j)\]
%*****************
\item $i \in O_P$ and $j \in E_H$. By definition $\oplus$, $i\oplus j=t+m\in O_H$ such that $t\in P(n)$ and $t\equiv i+j\pmod m$. By assumption
\begin{equation}\label{eq15}
A(i\oplus j)=A(t+m)\equiv(t+m)+\frac{m}{2} \equiv i+j+\frac{m}{2}\pmod m
\end{equation}
Also $A(j)=j$ and $A(i)=r\in O_P$  such that $r\equiv i+\frac{m}{2}\pmod m$.
By definition $\oplus$,
\begin{equation}\label{eq16}
A(i)\oplus A(j)=r\oplus j\equiv r+j\equiv i+\frac{m}{2}+j\pmod m
\end{equation}
Since $A(i)\oplus A(j),A(i\oplus j)\in H(n)$ by equations \ref{eq15}, \ref{eq16} and lemma \ref{lem2} $A(i\oplus j)=A(i)\oplus A(j)$.
%*****************
\item $i \in E_P$ and $j \in O_H$. The proof is similar to before case (c).
%By definition $\oplus$, $i\oplus j=t+m\in O_H$ such that $t\in P(n)$ and $t\equiv i+j\pmod m$. By assumption
%\begin{equation}\label{eq15}
%A(i\oplus j)=A(t+m)\equiv(t+m)+\frac{m}{2} \equiv i+j+\frac{m}{2}\pmod m
%\end{equation}
%Also $A(i)=i$ and $A(j)=r+m\in O_H$  such that $r\equiv j+\frac{m}{2}\pmod m$.
%By definition $\oplus$,
%\begin{equation}\label{eq16}
%A(i)\oplus A(j)=i\oplus (r+m)\equiv i+r\equiv i+j+\frac{m}{2}\pmod m
%\end{equation}
%Since $A(i)\oplus A(j),A(i\oplus j)\in H(n)$ by Equations \ref{eq15}, \ref{eq16} and Lemma \ref{lem2} $A(i\oplus j)=A(i)\oplus A(j)$.
\end{enumerate}
%%%%%%%%%%%%%%%%%%%%%%%%%%%%%%%%%%%%%
\item $i\in H(n)$ and $j\in P(n)$. It is the fourth case we consider in four steps as follows:
\begin{enumerate}
%*****************
\item $i \in O_H$ and $j \in O_P$. The proof is similar to (3)(a)  and so it is omitted.
%By definition $\oplus$, $i\oplus j=t+m\in E_H$ such that $t\in P(n)$ and $t\equiv i+j\pmod m$.
%By assumption
%\begin{equation}\label{eq17}
%A(i\oplus j)=A(t+m)=t+m\equiv  i+j\pmod m
%\end{equation}
%Also $A(i)=r_1+m\in O_H$ and $A(j)=r_2\in O_P$  such that
%\[r_1\equiv i+\frac{m}{2}\pmod m \hspace*{1cm} \text{and} \hspace*{1cm} r_2\equiv j+\frac{m}{2}\pmod m. \]
%By definition $\oplus$,
%\begin{equation}\label{eq18}
%A(i)\oplus A(j)=(r_1+m)\oplus r_2\equiv r_1+(r_2+m)\equiv i+j\pmod m
%\end{equation}
%Since $A(i)\oplus A(j),A(i\oplus j)\in H(n)$ by Equations \ref{eq17}, \ref{eq18} and Lemma \ref{lem2} $A(i\oplus j)=A(i)\oplus A(j)$.
%*****************
\item $i \in E_H$ and $j \in E_P$. The proof is similar to (3)(b)  and so it is omitted.
%It is clear $A(i)=i$ and $A(j)=j$. By definition $\oplus$, $i\oplus j=t+m\in E_H$ such that $t\in P(n)$ and $t\equiv i+j\pmod m$. Therefore
%\[A(i\oplus j)=A(t+m)=t+m=i\oplus j=A(i)\oplus A(j)\]
%*****************
\item $i \in O_H$ and $j \in E_P$. The proof is similar to (3)(c)  and so it is omitted.
%By definition $\oplus$, $i\oplus j=t+m\in O_H$ such that $t\in P(n)$ and $t\equiv i+j\pmod m$. By assumption
%\begin{equation}\label{eq15}
%A(i\oplus j)=A(t+m)\equiv(t+m)+\frac{m}{2} \equiv i+j+\frac{m}{2}\pmod m
%\end{equation}
%Also $A(j)=j$ and $A(i)=r+m\in O_H$  such that $r\equiv i+\frac{m}{2}\pmod m$.
%By definition $\oplus$,
%\begin{equation}\label{eq16}
%A(i)\oplus A(j)=(r+m)\oplus j\equiv r+j\equiv i+\frac{m}{2}+j\pmod m
%\end{equation}
%Since $A(i)\oplus A(j),A(i\oplus j)\in H(n)$ by Equations \ref{eq15}, \ref{eq16} and Lemma \ref{lem2} $A(i\oplus j)=A(i)\oplus A(j)$.
%*****************
\item $i \in E_H$ and $j \in O_P$. By definition $\oplus$, $i\oplus j=s+m\in O_H$ such that $s\in P(n)$ and $s\equiv i+j+\frac{m}{2}\pmod m$. By assumption
\begin{equation}\label{eq17}
A(i\oplus j)=A(s+m)\equiv(s+m)+\frac{m}{2} \equiv i+j\pmod m
\end{equation}
Also $A(i)=i$ and $A(j)=r\in O_P$  such that $r\equiv j+\frac{m}{2}\pmod m$.
By definition $\oplus$,
\begin{equation}\label{eq18}
A(i)\oplus A(j)=i\oplus r\equiv i+r+\frac{m}{2}\equiv i+j\pmod m
\end{equation}
Since $A(i)\oplus A(j),A(i\oplus j)\in H(n)$ by equations \ref{eq17}, \ref{eq18} and lemma \ref{lem2} $A(i\oplus j)=A(i)\oplus A(j)$.
\end{enumerate}
%%%%%%%%%%%%%%%%%%%%%%%%%%%%%%%%%%
\end{enumerate}
This completes the proof.
\end{proof}
%%%%%%%%%%%%%%%%%%%%%%%%%%%%%%%%%%
\begin{nota}\label{not2}
{\em
Assume Notation \ref{not1}. Set
$M=[O_P\times(O_H\cup E_H)]\bigcup[O_H\times(O_P\cup E_H)]\bigcup[E_H\times(O_P\cup O_H)]$. Define the map $\mathrm{gyr}:G_2(n)\times G_2(n)\longrightarrow \Aut{G_2(n),\oplus}$ as follows:
\[\mathrm{gyr}(a,b) = \gyr{a,b}{}=
\begin{cases}
A & :(a,b)\in M \\
I & \text{Otherwise}
\end{cases}\]
which $I, A\in \Aut{G_2(n),\oplus}$.  $I$ is identity automorphism and $A$ is a automorphism that is defined in lemma  \ref{lem1}. Obviously, the map $\mathrm{gyr}$ is well defined.
}
\end{nota}
%%%%%%%%%%%%%%%%%%%%%%%%%%%%%%%%%
%%%%%%%%%%%%%%%%%%%%%%%%%%%%%%%%%
%%%%%%%%%%%%%%%%%%%%%%%%%%%%%%%%%
\section{Main Result}
\noindent The aim of this section is to characterize finite gyrogroups of order $2^n$ constructed by the cyclic 2-group $\mathbb{Z}_{2^n}$, for $n\geq3$. We are now ready to state our main result.
%%%%%%%%%%%%%%%%%%%%%%%%%%%%%%%%%
\begin{thm}
Assume notations \ref{not1} and \ref{not2}. Then $(G_2(n),\oplus)$ is a gyrogroup.
\end{thm}

\begin{proof}
By Lemma \ref{lem1}, $\gyr{a,b}{}\in \Aut{G,\oplus}$. By definition of $\oplus$, $0\oplus i=i$ for all $i\in G_2(n)$ and so $0$ is the identity element of $G_2(n)$. If $\ominus x$ is the inverse of an arbitrary element $x\in G_2(n)$. we can be computed by the following formula:
\[\ominus x =
\begin{cases}
 -x & x\in P(n) \\
 -t+m & x=t+m\in H(n)
\end{cases}\]
in which $-x$ and $-t$ are inverse of $x$ and $t$ in $P(n)$, respectively.

Now we will prove the loop property. Assume $(a,b)$ be an arbitrary element of $G_2(n)\times G_2(n)$.
We will have four separate cases as follows:
\begin{enumerate}
\item $(a,b)\in P(n)\times P(n)$. In this case, $a\oplus b\in P(n)$. Clearly $(a\oplus b,b),(a,b)\not\in M$. By definition $\mathrm{gyr}$,
\[\gyr{a\oplus b,b}{}=I=\gyr{a,b}{}.\]
%****************
\item $(a,b)\in H(n)\times H(n)$. In this case, we have the following two subcases:
\begin{enumerate}
\item If $(a,b)\in (O_H\times O_H)\cup(E_H\times E_H)$, then $a\oplus b\in E_P$ and $(a,b),(a\oplus b,b)\not\in M$.
    Therefore \[\gyr{a\oplus b,b}{}=I=\gyr{a,b}{}.\]
\item If $(a,b)\in (O_H\times E_H)\cup(E_H\times O_H)$, then $a\oplus b\in O_P$ and $(a,b),(a\oplus b,b)\in M$.
    Therefore \[\gyr{a\oplus b,b}{}=A=\gyr{a,b}{}.\]
\end{enumerate}
%****************

\vspace{5cm}

\item $(a,b)\in P(n)\times H(n)$. In this case, $a\oplus b\in H$ and $a\oplus b\equiv a+b\pmod m$. Thus we have the following two subcases:
\begin{enumerate}
\item If $(a,b)\in (E_P\times O_H)\cup(E_P\times E_H)$, then $(a,b),(a\oplus b,b)\not\in M$.
    Therefore \[\gyr{a\oplus b,b}{}=I=\gyr{a,b}{}.\]
\item If $(a,b)\in (O_P\times O_H)\cup(O_P\times E_H)$, then $(a,b),(a\oplus b,b)\in M$.
    Therefore \[\gyr{a\oplus b,b}{}=A=\gyr{a,b}{}.\]
\end{enumerate}
%****************
\item $(a,b)\in H(n)\times P(n)$. In this case, $a\oplus b\in H$ and $a\oplus b\equiv a+b\pmod m$. Thus we have the following two subcases:
\begin{enumerate}
\item If $(a,b)\in (O_H\times E_P)\cup(E_H\times E_P)$, then $(a,b),(a\oplus b,b)\not\in M$.
    Therefore \[\gyr{a\oplus b,b}{}=I=\gyr{a,b}{}.\]
\item If $(a,b)\in (O_H\times O_P)\cup(E_H\times O_P)$, then $(a,b),(a\oplus b,b)\in M$.
    Therefore \[\gyr{a\oplus b,b}{}=A=\gyr{a,b}{}.\]
\end{enumerate}
%****************
\end{enumerate}
Therefore the loop property is valid.
%%%%%%%%%%%%%%%%%%%%%%%%%%%%%%%%%%
Finally, we investigate the left gyroassociative law. To do this, We have four separate cases as follows:
\begin{enumerate}
%%%%%%%%%%%%
\item $(a,b)\in P(n)\times P(n)$. In this case $a\oplus b\in P(n)$ and $gyr{a,b}{}=I$.
We consider two subcases as follows:
\begin{enumerate}
  \item $c\in P(n)$. Thus $b\oplus c\in P(n)$ and by definition $\oplus$:   \label{itm1}
\begin{eqnarray*}
(a\oplus b)\oplus \gyr{a,b}{c} & = & (a\oplus b)\oplus c \pmod m \notag\\
                              & \equiv & (a\oplus b)+c \pmod m \notag\\
							  & \equiv & a+b+c  \pmod m \\
							  & \equiv & a+(b\oplus c) \pmod m \notag\\
                              & \equiv & a\oplus(b\oplus c) \pmod m  \notag
\end{eqnarray*}
By Lemma \ref{lem2}, $a\oplus(b\oplus c)=(a\oplus b)\oplus \gyr{a,b}{c}$.
  \item $c\in H(n)$. Thus $b\oplus c\in H(n)$ and the proof of this case is similar to \ref{itm1}.
\end{enumerate}
%%%%%%%%%%%%%%%%%%%%%%%%%%%%%%%%%%%%%%5
\item $(a,b)\in H(n)\times H(n)$. By notation \ref{not2},
\[\gyr{a,b}{}=
\begin{cases}
A & :(a,b)\in N \\
I & \text{otherwise}
\end{cases}\]
in which $N=(O_H\times E_H)\cup(E_H\times O_H)\subseteq M$. We consider two subcases as follows:
\begin{enumerate}
%%%%%%%%%%%%%%%%%%%%
\item $(a,b)\not\in N$. In this subcase, $\gyr{a,b}{}=I$ and ($a,b\in O_H$ or $a,b\in E_H$).
If $a,b\in O_H$, then $a\oplus b\in E_P$. So we have the following cases:
\begin{enumerate}
  \item $c\in P(n)$. Thus $b\oplus c\in H(n)$ and the proof of this case is similar to \ref{itm1}.
  \item $c\in H(n)$. Thus $b\oplus c\in P(n)$ and the proof of this case is similar to \ref{itm1}.
  \end{enumerate}
If $a,b\in E_H$, then $a\oplus b\in E_P$. So we have the following cases:
\begin{enumerate}
\setcounter{enumiii}{2}
  \item $c\in O_P$. Thus $b\oplus c\in O_H$ and by definition $\oplus$:  \label{itm2}
  \begin{eqnarray*}
  (a\oplus b)\oplus \gyr{a,b}{c} & = & (a\oplus b)\oplus c \notag\\
									    & \equiv & (a\oplus b)+c \pmod m \notag\\
									    & \equiv & a+b+c \pmod m \\
									    & \equiv & a+(b+c+\frac{m}{2})+\frac{m}{2} \pmod m \\
									    & \equiv & a+(b\oplus c)+\frac{m}{2} \pmod m \\
									    & \equiv & a\oplus(b\oplus c) \pmod m
  \end{eqnarray*}
  By Lemma \ref{lem2}, $a\oplus(b\oplus c)=(a\oplus b)\oplus \gyr{a,b}{c}$.
  \item If $c\in E_P$, then $b\oplus c\in E_H$. The proof of this case is similar to \ref{itm1}.
  \item If $c\in O_H$, then $b\oplus c\in O_P$. The proof of this case is similar to \ref{itm2}.
  \item If $c\in E_H$, then $b\oplus c\in E_P$. The proof of this case is similar to \ref{itm1}.
  \end{enumerate}
%*****************************
\item $(a,b)\in N$. In this subcase $\gyr{a,b}{}=A$, also $(a\in O_H\  \& \ b\in E_H)$ or $(a\in E_H \  \& \  b\in O_H)$.
%%%%%%%%%%%%%%%%%%%%%%%%%%%%%%%%%%%%%%%%%%%%%%%%%%%%%%
    If $(a\in O_H\  \& \ b\in E_H)$, then $a\oplus b\in O_P$. So we have the following cases:
    \begin{enumerate}
    %****************
      \item $c\in O_P$. Thus $b\oplus c\in O_H$ and $gyr{a,b}{(c)}\in O_P$. By Lemma \ref{lem1} and definition $\oplus$: \label{itm3}
     \begin{eqnarray*}
     (a\oplus b)\oplus \gyr{a,b}{(c)} & \equiv & (a\oplus b)+\gyr{a,b}{(c)} \pmod m \notag\\
	           						    & \equiv & a+b+c+\frac{m}{2} \pmod m \notag\\
									    & \equiv & a+(b\oplus c) \pmod m \notag\\
									    & \equiv & a\oplus(b\oplus c) \pmod m
     \end{eqnarray*}
     By lemma \ref{lem2}, $a\oplus(b\oplus c)=(a\oplus b)\oplus \gyr{a,b}{(c)}$.
     %**********************
      \item $c\in E_P$. Thus $b\oplus c\in E_H$ and $\gyr{a,b}{(c)}=c$. The proof of this case is similar to \ref{itm1}.
     %*******************
     \item $c\in O_H$. Then $b\oplus c\in O_P$ and $\gyr{a,b}{(c)}\in O_H$. The proof of this case is similar to \ref{itm3}.
      %**************************
      \item $c\in E_H$. Then $b\oplus c\in E_P$ and $\gyr{a,b}{c}=c$. The proof of this case is similar to \ref{itm1}.
    \end{enumerate}
%%%%%%%%%%%%%%%%%%%%%%%%%%%%%%%%%%%%%%%%%%%%%%%%%%%%%%%
    If $(a\in E_H \  \& \  b\in O_H)$ then $a\oplus b\in O_P$. So we have the following cases:
    \begin{enumerate}
    \setcounter{enumiii}{4}
    %****************
      \item $c\in O_P$. Thus $b\oplus c\in E_H$ and $\gyr{a,b}{(c)}\in O_P$. By lemma \ref{lem1} and definition $\oplus$: \label{itm4}
     \begin{eqnarray*}
     (a\oplus b)\oplus \gyr{a,b}{(c)} & \equiv & (a\oplus b)+\gyr{a,b}{(c)} \pmod m \notag\\
	           						    & \equiv & a+b+\frac{m}{2}+c+\frac{m}{2} \pmod m \notag\\
	           						    & \equiv & a+b+c \pmod m \notag\\
									    & \equiv & a+(b\oplus c) \pmod m \notag\\
									    & \equiv & a\oplus(b\oplus c) \pmod m
     \end{eqnarray*}
     By lemma \ref{lem2}, $a\oplus(b\oplus c)=(a\oplus b)\oplus \gyr{a,b}{(c)}$.
     %**********************
      \item $c\in E_P$. Thus $b\oplus c\in O_H$ and $\gyr{a,b}{(c)}=c$. By definition $\oplus$: \label{itm5}
     \begin{eqnarray*}
     (a\oplus b)\oplus \gyr{a,b}{(c)} & = & (a\oplus b)\oplus c \pmod m \notag\\
									    & \equiv & (a\oplus b)+c \pmod m \notag\\
									    & \equiv & a+b+\frac{m}{2}+c \pmod m \notag\\
									    & \equiv & a+(b\oplus c)+\frac{m}{2} \pmod m \notag\\
									    & \equiv & a\oplus(b\oplus c) \pmod m
     \end{eqnarray*}
     By lemma \ref{lem2}, $a\oplus(b\oplus c)=(a\oplus b)\oplus \gyr{a,b}{(c)}$.
     %*******************
     \item $c\in O_H$. Then $b\oplus c\in E_P$ and $\gyr{a,b}{(c)}\in O_H$. The proof of this case is similar to \ref{itm4}.
      %**************************
      \item $c\in E_H$. Then $b\oplus c\in O_P$ and $\gyr{a,b}{(c)}=c$. The proof of this case is similar to \ref{itm5}.
    \end{enumerate}
\end{enumerate}
%%%%%%%%%%%%%%%%%%%%%%%%%%%%%%%%%%%%%%%%%%%%%
\item $(a,b)\in P(n)\times H(n)$. By notation \ref{not2}, it is clear
\[\gyr{a,b}{}=
\begin{cases}
A & :a\in O_P \\
I &  :a\not\in O_P
\end{cases}\]
We consider two subcases as follows:
\begin{enumerate}%[$(i)$]
\item $a\not\in O_P$. In this case $\gyr{a,b}{}=I$. If $b\in O_H$, then $a\oplus b\in O_H$. So we have the following cases:
    \begin{enumerate}
      \item $c\in P(n)$. Thus $b\oplus c\in H(n)$ and the proof of this case is similar to \ref{itm1}.
      \item $c\in H(n)$. Thus $b\oplus c\in P(n)$ and the proof of this case is similar to \ref{itm1}.
    \end{enumerate}
If $b\in E_H$, then $a\oplus b\in E_H$. So we have the following cases:
    \begin{enumerate}
    \setcounter{enumiii}{2}
      \item $c\in O_P$. Thus $b\oplus c\in O_H$ and by definition $\oplus$:\label{it02}
     \begin{eqnarray*}
     (a\oplus b)\oplus \gyr{a,b}{(c)} & = & (a\oplus b)\oplus c\\
                                        & \equiv & (a\oplus b)+ c+\frac{m}{2} \pmod m \\
									    & \equiv & a+b+c+\frac{m}{2} \pmod m \\
									    & \equiv & a+(b\oplus c) \pmod m \\
									    & \equiv & a\oplus(b\oplus c) \pmod m
     \end{eqnarray*}
     By Lemma \ref{lem2}, $a\oplus(b\oplus c)=(a\oplus b)\oplus \gyr{a,b}{(c)}$.
      \item $c\in E_P$. Thus $b\oplus c\in E_H$ and the proof of this case is similar to  \ref{itm1}.
      \item $c\in O_H$. Thus $b\oplus c\in O_P$ and the proof of this case is similar to \ref{it02}.
      \item $c\in E_H$. Thus $b\oplus c\in E_P$ and the proof of this case is similar to  \ref{itm1}.
     \end{enumerate}
%%%%%%%%%%%%%%%%%%%%
\item $a\in O_P$. In this case $\gyr{a,b}{}=A$. If $b\in O_H$, then $a\oplus b\in E_H$. So we have the following cases:
    \begin{enumerate}
      \item $c\in O_P$. Thus $b\oplus c\in E_H$ and $\gyr{a,b}{(c)}\in O_P$. By lemma \ref{lem1} and definition $\oplus$: \label{it03}
     \begin{eqnarray*}
     (a\oplus b)\oplus \gyr{a,b}{(c)}      & \equiv & (a\oplus b)+ \gyr{a,b}{(c)}+\frac{m}{2} \pmod m \\
									    & \equiv & a+b+c+\frac{m}{2}+\frac{m}{2} \pmod m \\
									    & \equiv & a+b+c \pmod m \\
									    & \equiv & a+(b\oplus c) \pmod m \\
									    & \equiv & a\oplus(b\oplus c) \pmod m
     \end{eqnarray*}
     By Lemma \ref{lem2}, $a\oplus(b\oplus c)=(a\oplus b)\oplus \gyr{a,b}{(c)}$.
      \item $c\in E_P$. Thus $b\oplus c\in O_H$ and $\gyr{a,b}{(c)}=c$. The proof of this case is similar to \ref{itm1}.
      \item $c\in O_H$. Thus $b\oplus c\in E_P$ and $\gyr{a,b}{(c)}\in O_H$. The proof of this case is similar to \ref{it03}.
      \item $c\in E_H$. Thus $b\oplus c\in O_P$ and $\gyr{a,b}{(c)}=c$. The proof of this case is similar to \ref{itm1}.
    \end{enumerate}
If $b\in E_H$, then $a\oplus b\in O_H$. So we have the following cases:
    \begin{enumerate}
    \setcounter{enumiii}{4}
      \item $c\in O_P$. Thus $b\oplus c\in O_H$ and $\gyr{a,b}{(c)}\in O_P$. The proof of this case is similar to \ref{itm3}.
      \item $c\in E_P$.Thus $b\oplus c\in E_H$ and $\gyr{a,b}{(c)}=c$. The proof of this case is similar to \ref{itm1}.
      \item $c\in O_H$. Thus $b\oplus c\in O_P$ and $\gyr{a,b}{(c)}\in O_H$. The proof of this case is similar to \ref{itm3}.
      \item $c\in E_H$. Thus $b\oplus c\in E_P$ and $\gyr{a,b}{(c)}=c$. and the proof of this case is similar to \ref{itm1}.
     \end{enumerate}
\end{enumerate}
%%%%%%%%%%%%%%%%%%%%%%%%%%%%%
%%%%%%%%%%%%%%%%%%%%%%%%%%%%%
\item $(a,b)\in H(n)\times P(n)$. By notation \ref{not2}, it is clear
\[\gyr{a,b}{}=
\begin{cases}
A & :b\in O_P \\
I &  :b\not\in O_P
\end{cases}\]
We consider two subcases as follows:
\begin{enumerate}
\item $b\not\in O_P$. In this case $\gyr{a,b}{}=I$.
If $a\in O_H$, then $a\oplus b\in O_H$. So we have the following cases:
    \begin{enumerate}
      \item $c\in P(n)$. Thus $b\oplus c\in P(n)$ and the proof of this case is similar to  \ref{itm1}.
      \item $c\in H(n)$. So $b\oplus c\in H(n)$ and the proof of this case is similar to  \ref{itm1}.
    \end{enumerate}
If $a\in E_H$, then $a\oplus b\in E_H$. So we have the following cases:
    \begin{enumerate}
    \setcounter{enumiii}{2}
      \item $c\in O_P$. Thus $b\oplus c\in O_P$ and by definition $\oplus$:\label{it05}
     \begin{eqnarray*}
     (a\oplus b)\oplus \gyr{a,b}{(c)} & = & (a\oplus b)\oplus c\\
                                        & \equiv & (a\oplus b)+ c+\frac{m}{2} \pmod m \\
									    & \equiv & a+b+c+\frac{m}{2} \pmod m \\
									    & \equiv & a+(b\oplus c)+\frac{m}{2} \pmod m \\
									    & \equiv & a\oplus(b\oplus c) \pmod m
     \end{eqnarray*}
     By Lemma \ref{lem2}, $a\oplus(b\oplus c)=(a\oplus b)\oplus \gyr{a,b}{(c)}$.
      \item $c\in E_P$.Thus $b\oplus c\in E_P$ and the proof of this case is similar to  \ref{itm1}.
      \item $c\in O_H$. Thus $b\oplus c\in O_H$ and the proof of this case is similar to \ref{it05}.
      \item $c\in E_H$. Thus $b\oplus c\in E_H$ and the proof of this case is similar to  \ref{itm1}.
     \end{enumerate}
%%%%%%%%%%%%%%%%%%%%
\item $b\in O_P$. In this case $\gyr{a,b}{}=A$.
If $a\in O_H$, then $a\oplus b\in E_H$. So we have the following cases:
    \begin{enumerate}
      \item $c\in O_P$. Thus $b\oplus c\in E_P$ and $\gyr{a,b}{(c)}\in O_P$. The proof of this case is similar to \ref{it03}.
      \item $c\in E_P$. Thus $b\oplus c\in O_P$ and $\gyr{a,b}{(c)}=c$. The proof of this case is similar to \ref{itm1}.
      \item $c\in O_H$. Thus $b\oplus c\in E_H$ and $\gyr{a,b}{(c)}\in O_H$. The proof of this case is similar to \ref{it03}.
      \item $c\in E_H$. Thus $b\oplus c\in O_H$ and $\gyr{a,b}{(c)}=c$. The proof of this case is similar to \ref{itm1}.
    \end{enumerate}
If $a\in E_H$, then $a\oplus b\in O_H$. So we have the following cases:
    \begin{enumerate}
    \setcounter{enumiii}{4}
      \item $c\in O_P$. Thus $b\oplus c\in E_P$ and $\gyr{a,b}{(c)}\in O_P$. The proof of this case is similar to \ref{itm4}.
      \item $c\in E_P$. Thus $b\oplus c\in O_P$ and $\gyr{a,b}{(c)}=c$. The proof of this case is similar to \ref{itm5}.
      \item $c\in O_H$. Thus $b\oplus c\in E_H$ and $\gyr{a,b}{(c)}\in O_H$. The proof of this case is similar to \ref{itm4}.
      \item $c\in E_H$. Thus $b\oplus c\in O_H$ and $\gyr{a,b}{(c)}=c$. and the proof of this case is similar to \ref{itm5}.
     \end{enumerate}
\end{enumerate}
%%%%%%%%%%%%%%%%%%%
%%%%%%%%%%%%%%%%%%%
\end{enumerate}
By the above mentions, $(G_2(n),\oplus)$ is a grogroup and this completes the proof.
\end{proof}
%%%%%%%%%%%%%%%%%%%%%%%%%%%%%%%%%%%
Ashrafi et al. in Theorem 2 \cite{12} proved the gyrogroup $G_1(n)$ is a non-gyrocommutative, but in the following theorem we show the gyrogroup $G_2(n)$ is a gyrocommutative. Therfore the gyrogroups $G_1(n)$ and $G_2(n)$ are not isomorphic.
%%%%%%%%%%%%%%%%%%%%%%%%%%%%%%%%%
\begin{thm}
The gyrogroup $(G_2(n),\oplus)$ is a gyrocommutative.
\end{thm}

\begin{proof}
Suppose $(a,b)$ is an arbitrary member of $G_2(n)\times G_2(n)$. We consider two cases as follows:
\begin{enumerate}
  \item $(a,b)\not\in M$. In this case $\gyr{a,b}=I$. By definition $\oplus$,
      \[\gyr{a,b}{(b\oplus a)}=b\oplus a \equiv a+b \equiv a\oplus b \pmod m \]
      By lemma \ref{lem2}, $\gyr{a,b}{(b\oplus a)}=a\oplus b$.
  \item $(a,b)\in M$. In this case $\gyr{a,b}{}=A$. If $b\oplus a$ is an even number, then $(a,b)\in S\cup S^{-1}$ in which $S=O_P\times O_H$. By lemma \ref{lem1} and definition $\oplus$,
      \[\gyr{a,b}{(b\oplus a)}=b\oplus a \equiv a+b \equiv a\oplus b \pmod m \]
      By lemma \ref{lem2}, $\gyr{a,b}{(b\oplus a)}=a\oplus b$. If $b\oplus a$ is an odd number,
      then $(a,b)\in T\cup T^{-1}$ in which $T=E_H\times(O_P\cup O_H)$. Now, we consider two subcases as follows:
      \begin{enumerate}[$(i)$]
        \item $(a,b)\in T$. So $(b,a)\in T^{-1}$. By lemma \ref{lem1} and definition $\oplus$,
            \[\gyr{a,b}{(b\oplus a)}\equiv b\oplus a+\frac{m}{2} \equiv a+b+\frac{m}{2} \equiv a\oplus b \pmod m \]
            By lemma \ref{lem2}, $\gyr{a,b}{(b\oplus a)}=a\oplus b$.
        \item $(a,b)\in T^{-1}$. So $(b,a)\in T$. By lemma \ref{lem1} and definition $\oplus$,
            \[\gyr{a,b}{(b\oplus a)}\equiv b\oplus a+\frac{m}{2} \equiv a+b+\frac{m}{2}+\frac{m}{2}\equiv a+b \equiv a\oplus b \pmod m \]
            By lemma \ref{lem2}, $\gyr{a,b}{(b\oplus a)}=a\oplus b$.
      \end{enumerate}
\end{enumerate}
This completes the proof of gyrocommutativity.
\end{proof}
%%%%%%%%%%%%%%%%%%%%%%%%%%%%%%%%%%%%%%%%%%%%%%%%%%%%%%%%%%
%%%%%%%%%%%%%%%%%%%%%%%%%%%%%%%%%%%
\begin{exa}\label{ex1}
In this example, we investigate the gyrogroup $G_2(3)$ of order $8$ constructed by cyclic group $\mathbb{Z}_8$. By definition, $G_2(3)=\{0,1,2,3,4,5,6,7\}$ and the binary operation $\oplus$ is defined as follows:
\[i\oplus j =
\begin{cases}
t & :(i,j)\in (P(3)\times P(3))\cup\Bigl[(H(3)\times H(3))-(E_H\times O_H)\Bigr]\\[2mm]
t+4 & :(i,j)\in (P(3)\times H(3))\cup\Bigl[(H(3)\times P(3))-(E_H\times O_P)\Bigr]\\[1mm]
s & :(i,j)\in E_H\times O_H \\[1mm]
s+4 & :(i,j)\in E_H\times O_P
\end{cases}\]
in which $P(3)=\{0,1,2,3\}$, $H(3)=\{4,5,6,7\}$ and  $t,s\in P(3)$. Also
\[\left\{\begin{split}
t & \equiv i+j \pmod 4\\
s & \equiv i+j+2 \pmod 4
\end{split}
\right..\]
Therefore the addition table and the gyration table for $G_2(3)$ are presented in Table \ref{tab1} in which $A$ is
the unique non-identity gyroautomorphism of $G_2(3)$ given by $A = (1, 3)(5, 7)$.
\begin{table}[htp]
\centering \caption{\label{tab1} \footnotesize }
\renewcommand*{\arraystretch}{1}
\subfloat[][The Cayley table]{
     \begin{tabular}{ c | c  c  c  c  c  c  c  c }
            
             $\oplus$ & $0$ & $1$ & $2$ & $3$ & $4$ & $5$ & $6$ & $7$\\
              \hline
                   $0$ & $0$ & $1$ & $2$ & $3$ & $4$ & $5$ & $6$ & $7$\\
            
                   $1$ & $1$ & $2$ & $3$ & $0$ & $5$ & $6$ & $7$ & $4$\\
            
                   $2$ & $2$ & $3$ & $0$ & $1$ & $6$ & $7$ & $4$ & $5$\\
            
                   $3$ & $3$ & $0$ & $1$ & $2$ & $7$ & $4$ & $5$ & $6$\\
           
                   $4$ & $4$ & $7$ & $6$ & $5$ & $0$ & $3$ & $2$ & $1$\\
            
                   $5$ & $5$ & $6$ & $7$ & $4$ & $1$ & $2$ & $3$ & $0$\\
            
                   $6$ & $6$ & $5$ & $4$ & $7$ & $2$ & $1$ & $0$ & $3$\\
            
                   $7$ & $7$ & $4$ & $5$ & $6$ & $3$ & $0$ & $1$ & $2$\\
            
     \end{tabular}
}
\hspace*{3cm}
\subfloat[][The gyration table]{
     \begin{tabular}{ c | c  c  c  c  c  c  c  c }
            
             $\mathrm{gyr}$ & $0$ & $1$ & $2$ & $3$ & $4$ & $5$ & $6$ & $7$\\
              \hline
                $0$ & $I$ & $I$ & $I$ & $I$ & $I$ & $I$ & $I$ & $I$\\
            
                $1$ & $I$ & $I$ & $I$ & $I$ & $A$ & $A$ & $A$ & $A$\\
            
                $2$ & $I$ & $I$ & $I$ & $I$ & $I$ & $I$ & $I$ & $I$\\
            
                $3$ & $I$ & $I$ & $I$ & $I$ & $A$ & $A$ & $A$ & $A$\\
           
                $4$ & $I$ & $A$ & $I$ & $A$ & $I$ & $A$ & $I$ & $A$\\
            
                $5$ & $I$ & $A$ & $I$ & $A$ & $A$ & $I$ & $A$ & $I$\\
            
                $6$ & $I$ & $A$ & $I$ & $A$ & $I$ & $A$ & $I$ & $A$\\
            
                $7$ & $I$ & $A$ & $I$ & $A$ & $A$ & $I$ & $A$ & $I$\\
            
     \end{tabular}
}
\end{table}
\end{exa}
%%%%%%%%%%%%%%%%%%%%%%%%%%%%%%%%%%
\begin{exa}
In this example, we investigate the gyrogroup $G_2(4)$ of order $16$ constructed by cyclic group $\mathbb{Z}_{16}$. 
By definition, $G_2(4)=\{0,1,2,3,4,5,6,7,8,9,10,11,12,13,14,15\}$ and the binary operation $\oplus$ is defined as follows:
\[i\oplus j =
\begin{cases}
t & :(i,j)\in (P(4)\times P(4))\cup\Bigl[(H(4)\times H(4))-(E_H\times O_H)\Bigr]\\[2mm]
t+8 & :(i,j)\in (P(4)\times H(4))\cup\Bigl[(H(4)\times P(4))-(E_H\times O_P)\Bigr]\\[1mm]
s & :(i,j)\in E_H\times O_H \\[1mm]
s+8 & :(i,j)\in E_H\times O_P
\end{cases}\]
which $P(4)=\{0,1,2,3,4,5,6,7\}$, $H(4)=\{8,9,10,11,12,13,14,15\}$ and  $t,s\in P(4)$. Also
\[\left\{\begin{split}
t & \equiv i+j \pmod 8\\
s & \equiv i+j+4 \pmod 8
\end{split}
\right..\]
Therefore the addition table and the gyration table for $G_2(4)$ are presented in tables \ref{tab2} and \ref{tab3}, respectively, 
in which $A$ is the unique non-identity gyroautomorphism of $G_2(4)$ given by $A=(1, 5)(3, 7)(9, 13)(11, 15)$.
\begin{table}[htp]
\centering \caption{The Cayley table of $G_2(4)$ \label{tab2} \footnotesize }
 \footnotesize
 \begin{tabular}{ c | c  c  c  c  c  c  c  c  c  c  c  c  c  c  c  c}
 
 $\oplus$   & $0$ & $1$ & $2$ & $3$ & $4$ & $5$ & $6$ & $7$ & $8$ & $9$ & $10$ & $11$ & $12$ & $13$ & $14$ & $15$ \\ \hline
  $0$ & $0$ & $1$ & $2$ & $3$ & $4$ & $5$ & $6$ & $7$ & $8$ & $9$ & $10$ & $11$ & $12$ & $13$ & $14$ & $15$ \\
  $1$ & $1$ & $2$ & $3$ & $4$ & $5$ & $6$ & $7$ & $0$ & $9$ & $10$ & $11$ & $12$ & $13$ & $14$ & $15$ & $8$ \\
  $2$ & $2$ & $3$ & $4$ & $5$ & $6$ & $7$ & $0$ & $1$ & $10$ & $11$ & $12$ & $13$ & $14$ & $15$ & $8$ & $9$ \\
  $3$ & $3$ & $4$ & $5$ & $6$ & $7$ & $0$ & $1$ & $2$ & $11$ & $12$ & $13$ & $14$ & $15$ & $8$ & $9$ & $10$ \\
  $4$ & $4$ & $5$ & $6$ & $7$ & $0$ & $1$ & $2$ & $3$ & $12$ & $13$ & $14$ & $15$ & $8$ & $9$ & $10$ & $11$ \\
  $5$ & $5$ & $6$ & $7$ & $0$ & $1$ & $2$ & $3$ & $4$ & $13$ & $14$ & $15$ & $8$ & $9$ & $10$ & $11$ & $12$ \\
  $6$ & $6$ & $7$ & $0$ & $1$ & $2$ & $3$ & $4$ & $5$ & $14$ & $15$ & $8$ & $9$ & $10$ & $11$ & $12$ & $13$ \\
  $7$ & $7$ & $0$ & $1$ & $2$ & $3$ & $4$ & $5$ & $6$ & $15$ & $8$ & $9$ & $10$ & $11$ & $12$ & $13$ & $14$ \\
  $8$ & $8$ & $13$ & $10$ & $15$ & $12$ & $9$ & $14$ & $11$ & $0$ & $5$ & $2$ & $7$ & $4$ & $1$ & $6$ & $3$ \\
  $9$ & $9$ & $10$ & $11$ & $12$ & $13$ & $14$ & $15$ & $8$ & $1$ & $2$ & $3$ & $4$ & $5$ & $6$ & $7$ & $0$ \\
  $10$ & $10$ & $15$ & $12$ & $9$ & $14$ & $11$ & $8$ & $13$ & $2$ & $7$ & $4$ & $1$ & $6$ & $3$ & $0$ & $5$\\
  $11$ & $11$ & $12$ & $13$ & $14$ & $15$ & $8$ & $9$ & $10$ & $3$ & $4$ & $5$ & $6$ & $7$ & $0$ & $1$ & $2$\\
  $12$ & $12$ & $9$ & $14$ & $11$ & $8$ & $13$ & $10$ & $15$ & $4$ & $1$ & $6$ & $3$ & $0$ & $5$ & $2$ & $7$\\
  $13$ & $13$ & $14$ & $15$ & $8$ & $9$ & $10$ & $11$ & $12$ & $5$ & $6$ & $7$ & $0$ & $1$ & $2$ & $3$ & $4$\\
  $14$ & $14$ & $11$ & $8$ & $13$ & $10$ & $15$ & $12$ & $9$ & $6$ & $3$ & $0$ & $5$ & $2$ & $7$ & $4$ & $1$\\
  $15$ & $15$ & $8$ & $9$ & $10$ & $11$ & $12$ & $13$ & $14$ & $7$ & $0$ & $1$ & $2$ & $3$ & $4$ & $5$ & $6$\\
 
 \end{tabular}
\end{table}

\begin{table}[htp]
\centering \caption{The gyration table of $G_2(4)$ \label{tab3} \footnotesize }
 \footnotesize
 \begin{tabular}{ c | c  c  c  c  c  c  c  c  c  c  c  c  c  c  c  c}

 $\mathrm{gyr}$   & $0$ & $1$ & $2$ & $3$ & $4$ & $5$ & $6$ & $7$ & $8$ & $9$ & $10$ & $11$ & $12$ & $13$ & $14$ & $15$ \\ \hline
  $0$ & $I$ & $I$ & $I$ & $I$ & $I$ & $I$ & $I$ & $I$ & $I$ & $I$ & $I$ & $I$ & $I$ & $I$ & $I$ & $I$ \\
  $1$ & $I$ & $I$ & $I$ & $I$ & $I$ & $I$ & $I$ & $I$ & $A$ & $A$ & $A$ & $A$ & $A$ & $A$ & $A$ & $A$ \\
  $2$ & $I$ & $I$ & $I$ & $I$ & $I$ & $I$ & $I$ & $I$ & $I$ & $I$ & $I$ & $I$ & $I$ & $I$ & $I$ & $I$ \\
  $3$ & $I$ & $I$ & $I$ & $I$ & $I$ & $I$ & $I$ & $I$ & $A$ & $A$ & $A$ & $A$ & $A$ & $A$ & $A$ & $A$ \\
  $4$ & $I$ & $I$ & $I$ & $I$ & $I$ & $I$ & $I$ & $I$ & $I$ & $I$ & $I$ & $I$ & $I$ & $I$ & $I$ & $I$ \\
  $5$ & $I$ & $I$ & $I$ & $I$ & $I$ & $I$ & $I$ & $I$ & $A$ & $A$ & $A$ & $A$ & $A$ & $A$ & $A$ & $A$ \\
  $6$ & $I$ & $I$ & $I$ & $I$ & $I$ & $I$ & $I$ & $I$ & $I$ & $I$ & $I$ & $I$ & $I$ & $I$ & $I$ & $I$ \\
  $7$ & $I$ & $I$ & $I$ & $I$ & $I$ & $I$ & $I$ & $I$ & $A$ & $A$ & $A$ & $A$ & $A$ & $A$ & $A$ & $A$ \\
  $8$ & $I$ & $A$ & $I$ & $A$ & $I$ & $A$ & $I$ & $A$ & $I$ & $A$ & $I$ & $A$ & $I$ & $A$ & $I$ & $A$ \\
  $9$ & $I$ & $A$ & $I$ & $A$ & $I$ & $A$ & $I$ & $A$ & $A$ & $I$ & $A$ & $I$ & $A$ & $I$ & $A$ & $I$ \\
  $10$ & $I$ & $A$ & $I$ & $A$ & $I$ & $A$ & $I$ & $A$ & $I$ & $A$ & $I$ & $A$ & $I$ & $A$ & $I$ & $A$ \\
  $11$ & $I$ & $A$ & $I$ & $A$ & $I$ & $A$ & $I$ & $A$ & $A$ & $I$ & $A$ & $I$ & $A$ & $I$ & $A$ & $I$ \\
  $12$ & $I$ & $A$ & $I$ & $A$ & $I$ & $A$ & $I$ & $A$ & $I$ & $A$ & $I$ & $A$ & $I$ & $A$ & $I$ & $A$ \\
  $13$ & $I$ & $A$ & $I$ & $A$ & $I$ & $A$ & $I$ & $A$ & $A$ & $I$ & $A$ & $I$ & $A$ & $I$ & $A$ & $I$ \\
  $14$ & $I$ & $A$ & $I$ & $A$ & $I$ & $A$ & $I$ & $A$ & $I$ & $A$ & $I$ & $A$ & $I$ & $A$ & $I$ & $A$ \\
  $15$ & $I$ & $A$ & $I$ & $A$ & $I$ & $A$ & $I$ & $A$ & $A$ & $I$ & $A$ & $I$ & $A$ & $I$ & $A$ & $I$ \\
 
 \end{tabular}
\end{table}
\end{exa}
%%%%%%%%%%%%%%%%%%%%%%%%%%%%%%%%%%%%%%%%%%%%%%%
\begin{exa}
By GAP, we obtain all proper subgyrogroups of $G_2(3)$ and $G_2(4)$ such that are subgroup. 
%The lattice of two gyrogroups is shown in Figure~\ref{fig1}. 
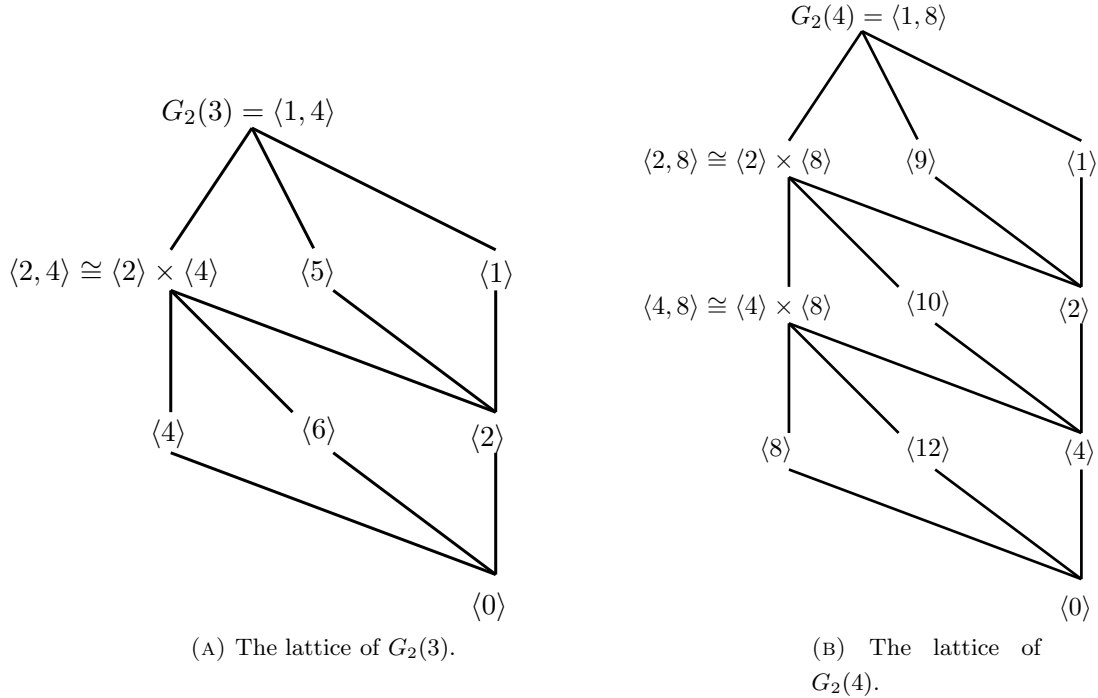
\begin{figure}[ht!]
\centering
\caption{The lattice of $G(3)$ and $G(4)$.}\label{fig1}
\subfloat[The lattice of $G_2(3)$.]{
\scalebox{1} % Change this value to rescale the drawing.
{
\begin{pspicture}(0,-3.4729688)(5.6228123,3.4729688)
\psline[linewidth=0.04cm](0.7809375,1.4545312)(1.8609375,3.0745313)
\psline[linewidth=0.04cm](0.7809375,0.91453123)(0.7809375,-0.7054688)
\psline[linewidth=0.04cm](5.1009374,0.91453123)(5.1009374,-0.7054688)
\psline[linewidth=0.04cm](5.1009374,-1.2454687)(5.1009374,-2.8654687)
\psline[linewidth=0.04cm](2.9409375,-1.2454687)(5.1009374,-2.8654687)
\psline[linewidth=0.04cm](0.7809375,-1.2454687)(5.1009374,-2.8654687)
\psline[linewidth=0.04cm](0.7809375,0.91453123)(5.1009374,-0.7054688)
\psline[linewidth=0.04cm](2.9409375,0.91453123)(5.1009374,-0.7054688)
\psline[linewidth=0.04cm](0.7809375,0.91453123)(2.4009376,-0.7054688)
\psline[linewidth=0.04cm](5.1009374,1.4545312)(1.8609375,3.0745313)
\psline[linewidth=0.04cm](2.6809375,1.4745313)(1.8609375,3.0745313)
\usefont{T1}{ptm}{m}{n}
\rput(5.012344,-3.2954688){$\langle 0\rangle$}
\usefont{T1}{ptm}{m}{n}
\rput(4.992344,-1.0554688){$\langle 2\rangle$}
\usefont{T1}{ptm}{m}{n}
\rput(2.7523437,-0.9554688){$\langle 6\rangle$}
\usefont{T1}{ptm}{m}{n}
\rput(0.7523438,-0.99546874){$\langle 4\rangle$}
\usefont{T1}{ptm}{m}{n}
\rput(0.023437,1.1445312){$\langle 2,4\rangle\cong\langle 2\rangle\times\langle 4\rangle$}
\usefont{T1}{ptm}{m}{n}
\rput(2.7323437,1.1445312){$\langle 5\rangle$}
\usefont{T1}{ptm}{m}{n}
\rput(5.112344,1.1245313){$\langle 1\rangle$}
\usefont{T1}{ptm}{m}{n}
\rput(1.8223437,3.2845314){$G_2(3)=\langle 1,4\rangle$}
\end{pspicture}
}
}
\hspace{2cm}
%%%%%%%%%%%%%%%%%%%%%%%
\subfloat[The lattice of $G_2(4)$.]{
\scalebox{.9} % Change this value to rescale the drawing.
{
\begin{pspicture}(0,-4.552969)(5.6228123,4.552969)
\psline[linewidth=0.04cm](0.7809375,2.5345314)(1.8609375,4.1545315)
\psline[linewidth=0.04cm](0.7809375,1.9945313)(0.7809375,0.37453124)
\psline[linewidth=0.04cm](5.1009374,1.9945313)(5.1009374,0.37453124)
\psline[linewidth=0.04cm](5.1009374,-2.3254688)(5.1009374,-3.9454687)
\psline[linewidth=0.04cm](2.9409375,-2.3254688)(5.1009374,-3.9454687)
\psline[linewidth=0.04cm](0.7809375,-2.3254688)(5.1009374,-3.9454687)
\psline[linewidth=0.04cm](0.7809375,1.9945313)(5.1009374,0.37453124)
\psline[linewidth=0.04cm](2.9409375,1.9945313)(5.1009374,0.37453124)
\psline[linewidth=0.04cm](0.7809375,1.9945313)(2.4009376,0.37453124)
\psline[linewidth=0.04cm](5.1009374,2.5345314)(1.8609375,4.1545315)
\psline[linewidth=0.04cm](2.6809375,2.5545313)(1.8609375,4.1545315)
\usefont{T1}{ptm}{m}{n}
\rput(5.012344,-4.3754687){$\langle 0\rangle$}
\usefont{T1}{ptm}{m}{n}
\rput(4.992344,0.02453125){$\langle 2\rangle$}
\usefont{T1}{ptm}{m}{n}
\rput(2.8423438,0.12453125){$\langle 10\rangle$}
\usefont{T1}{ptm}{m}{n}
\rput(5.072344,-2.0754688){$\langle 4\rangle$}
\usefont{T1}{ptm}{m}{n}
\rput(0.023437,2.2245312){$\langle 2,8\rangle\cong\langle 2\rangle\times\langle 8\rangle$}
\usefont{T1}{ptm}{m}{n}
\rput(2.7323437,2.2245312){$\langle 9\rangle$}
\usefont{T1}{ptm}{m}{n}
\rput(5.112344,2.2045312){$\langle 1\rangle$}
\usefont{T1}{ptm}{m}{n}
\rput(1.9723437,4.364531){$G_2(4)=\langle 1,8\rangle$}
\psline[linewidth=0.04cm](5.1009374,-0.16546875)(5.1009374,-1.7854687)
\psline[linewidth=0.04cm](0.7809375,-0.16546875)(5.1009374,-1.7854687)
\psline[linewidth=0.04cm](2.9409375,-0.16546875)(5.1009374,-1.7854687)
\psline[linewidth=0.04cm](0.7809375,-0.16546875)(2.4009376,-1.7854687)
\usefont{T1}{ptm}{m}{n}
\rput(2.8423438,-2.0354688){$\langle 12\rangle$}
\usefont{T1}{ptm}{m}{n}
\rput(0.59234375,-2.0354688){$\langle 8\rangle$}
\usefont{T1}{ptm}{m}{n}
\rput(0.023437,0.06453125){$\langle 4,8\rangle\cong\langle 4\rangle\times\langle 8\rangle$}
\psline[linewidth=0.04cm](0.7809375,-0.16546875)(0.7809375,-1.7854687)
\end{pspicture}
}
}
\end{figure}
%%%%%%%%%%%%%%%%%%%%%%%%%%%%%%%%%%%%%%%%%%%%%5
%Also, we guess that the lattice of  $G_2(n)$ is as Figure \ref{fig2}.

%%%%%%%%%%%%%%%%%%%%%%%%%%%%%%%%%%%%%%%%%%%
\end{exa}
%%%%%%%%%%%%%%%%%%%%%%%%%%%%%%%%%%%%%%%%%%%%%
\begin{remark}\label{Inn}
By using the gyrator identity, $\gyr{a,b}{c} = \ominus(a\oplus  b)\oplus(a\oplus (b\oplus c))$, one can simply calculate the group of all gyroautomorphisms of $G_2(4)$ is isomorphic $\mathbb{Z}_2$.
\end{remark}
\begin{proof}
Since the order of the unique non-identity gyroautomorphism of $G_2(4)$ is 2, it can forms a group with identity under the composition of permutations. This group will be isomorphic to $\mathbb{Z}_2$.
\end{proof}
\begin{prop}
By using construction of gyrosemidirect product of gyrogroups in \cite{8} one can calculate the gyroholomorph of $G_2(4)$ which is have the following structure,
$$\mathbb{Z}_2\times (\mathbb{Z}_8\rtimes\mathbb{Z}_2).$$
\end{prop}
\begin{thm}
Let $G_2(n)$ be the 2-gurogroup of order $2^n$. $H$ is subgyrogroup of  $G_2(n)$ if and only if it has one of the following forms:
\begin{enumerate}
\item $H=\langle 2^s\rangle$ is a subgroup of $P(n)$ such that $0\leq s\leq n-1$.
\item $H=\langle 2^s, m\rangle$ such that $0\leq s\leq n-1$ and $m=2^{n-1}$. All $H$'s are subgroups of $G_2(n)$ except 
$H=\langle 1, m\rangle=G_2(n)$
\item $H=\langle m+2^s\rangle$ is a subgroup of $G_2(n)$ such that $0\leq s\leq n-2$ and $m=2^{n-1}$.
\end{enumerate}
\end{thm}
%%%%%%%%%%%%%%%%%%%%%%
\begin{proof}
We know $G=P(n)\cup H(n)$ such that $P(n)\cap H(n)=\emptyset$. 
Since the restriction $\oplus$ of to $P(n)$ is the group addition, $P(n)\cong \mathbb{Z}_m$ that $m=2^{n-1}$. 
Suppose $H$ is a subgyrogroup of $G(n)$. 
If $H\subseteq P(n)$, then $H$ will be a subgroup of $P(n)$ and $H=\langle 2^s \rangle$ such that $0\leq s\leq n-1$, as desired. 
But $H\nsubseteq P(n)$, then $H=H_1\cup H_2$ in which $H_1\subseteq P(n)$ and $H_2\subseteq H(n)$. It is easy to see $H_1$ is a subgroup of $P(n)$. We consider two cases as follows: 
\begin{enumerate}
\item $m\in H_2$. In this case, the map $L_m\colon H_1 \to H_2$ defined by $L_m(x)=m\oplus x$  is bijective. 
Therefore $H_2=m\oplus H_1$. Since $H_1\leq P(n)\cong \mathbb{Z}_m$, $H=H_1\cup H_2=H_1\cup (m\oplus H_1)=\langle 2^s,m\rangle$ 
such that $0\leq s\leq n-1$. For $s=0$ and $s=n-1$, $\langle 1,m\rangle=G(n)$ and $\langle 0,m\rangle=\langle m\rangle$, respectively. 
If $s\neq 0$, then by the definition of gyroautomorphism of $G_2(n)$, all gyroautomorphisms of $H$ are the identity automorphism. Therefore $H$'s are subgroups of $G_2(n)$ and $H\cong\langle 2^s\rangle\times\langle m\rangle$.
\item $m\not\in H_2$. There exists an integer $i\in \mathbb{Z}_m$ such that $i\neq 0$ and $i$ is the smallest number in which $m+i \in H_2$. 
Also, the map $L_{m+i}\colon H_1 \to H_2$ defined by $L_{m+i}(x)=(m+i)\oplus x$  is bijective. Therefore $H_2=(m+i)\oplus H_1$. 
By definition of $\oplus$, $(m+i)\oplus(m+i)=2i\in H_1$. Since $i$ is the smallest number that has been chosen 
and $H_1\leq P(n)\cong \mathbb{Z}_m$, then $H_1=\langle 2i\rangle$ and 
$H=H_1\cup H_2=H_1\cup \Big((m+i)\oplus H_1\Big)=\langle 2i,m+i\rangle=\langle m+i\rangle$ such that $i=2^s$ and $0\leq s\leq n-2$.
\end{enumerate}
It is easy, $\langle 2^s\rangle\leq \langle 2^s,m\rangle$, $\langle m+2^s\rangle\leq \langle 2^s,m\rangle$ and  
$\langle m\rangle\leq \langle 2^{n-2},m\rangle\leq \langle 2^{n-3},m\rangle\leq\cdots \leq\langle 2,m\rangle\leq \langle 1,m\rangle$ thus the subgyrogroup lattices of $G(n)$ are depicted in Figure \ref{fig2} 
%%%%%%%%%%%%%%%%%%%%%%%%%%%%%%%%%%%%%%%%%%%%%%%%%%
\begin{figure}[ht!]
\centering % \footnotesize
\caption{The lattice of $G_2(n)$.}\label{fig2}
~\vspace*{.5cm}\\
\scalebox{.9} % Change this value to rescale the drawing.
{
\begin{pspicture}(0,-5.902969)(5.9628124,5.902969)
\psline[linewidth=0.04cm](0.7809375,3.8845313)(1.8609375,5.5045314)
\psline[linewidth=0.04cm](0.7809375,0.64453125)(0.7809375,-0.97546875)
\psline[linewidth=0.04cm](5.1009374,3.3445313)(5.1009374,1.7245313)
\psline[linewidth=0.04cm](5.1009374,-3.6754687)(5.1009374,-5.295469)
\psline[linewidth=0.04cm](2.9409375,-3.6754687)(5.1009374,-5.295469)
\psline[linewidth=0.04cm](0.7809375,-3.6754687)(5.1009374,-5.295469)
\psline[linewidth=0.04cm](0.7809375,0.64453125)(5.1009374,-0.97546875)
\psline[linewidth=0.04cm](2.9409375,3.3445313)(5.1009374,1.7245313)
\psline[linewidth=0.04cm](0.7809375,0.64453125)(2.4009376,-0.97546875)
\psline[linewidth=0.04cm](5.1009374,3.8845313)(1.8609375,5.5045314)
\psline[linewidth=0.04cm](2.6809375,3.9045312)(1.8609375,5.5045314)
\usefont{T1}{ptm}{m}{n}
\rput(5.012344,-5.7254686){$\langle 0\rangle$}
\usefont{T1}{ptm}{m}{n}
\rput(5.182344,-1.3254688){$\langle 2^{n-3}\rangle$}
\usefont{T1}{ptm}{m}{n}
\rput(3.1823437,-1.2254688){$\langle m+2^{n-3}\rangle$}
\usefont{T1}{ptm}{m}{n}
\rput(5.262344,-3.4254687){$\langle 2^{n-2}\rangle$}
\usefont{T1}{ptm}{m}{n}
\rput(.05234374,3.5745313){$\langle 2,m\rangle\cong\langle 2\rangle\times\langle m\rangle$}
\usefont{T1}{ptm}{m}{n}
\rput(2.9723437,3.5745313){$\langle m+1\rangle$}
\usefont{T1}{ptm}{m}{n}
\rput(5.112344,3.5545313){$\langle 1\rangle$}
\usefont{T1}{ptm}{m}{n}
\rput(1.8223437,5.7145314){$G_2(n)=\langle 1,m\rangle$}
\psline[linewidth=0.04cm](5.1009374,-1.5154687)(5.1009374,-3.1354687)
\psline[linewidth=0.04cm](0.7809375,-1.5154687)(5.1009374,-3.1354687)
\psline[linewidth=0.04cm](2.9409375,-1.5154687)(5.1009374,-3.1354687)
\psline[linewidth=0.04cm](0.7809375,-1.5154687)(2.4009376,-3.1354687)
\usefont{T1}{ptm}{m}{n}
\rput(2.6423438,-3.3854687){$\langle m+2^{n-2}\rangle$}
\usefont{T1}{ptm}{m}{n}
\rput(0.64234376,-3.3854687){$\langle m\rangle$}
\usefont{T1}{ptm}{m}{n}
\rput(-.50234374,-1.2854687){$\langle 2^{n-2},m\rangle\cong\langle 2^{n-2}\rangle\times\langle m \rangle$}
\psline[linewidth=0.04cm](0.7809375,-1.5154687)(0.7809375,-3.1354687)
\usefont{T1}{ptm}{m}{n}
\rput(-.50234374,0.87453127){$\langle 2^{n-3},m\rangle\cong\langle 2^{n-3}\rangle\times\langle m\rangle$}
\psline[linewidth=0.04cm](5.1009374,1.7245313)(0.7809375,3.3445313)
\psline[linewidth=0.04cm,linestyle=dotted,dotsep=0.16cm](0.7809375,3.3445313)(0.7809375,1.1845312)
\psline[linewidth=0.04cm,linestyle=dotted,dotsep=0.16cm](5.1009374,1.1845312)(5.1009374,-0.97546875)
\usefont{T1}{ptm}{m}{n}
\rput(5.112344,1.3945312){$\langle 2\rangle$}
\end{pspicture}
}
\end{figure}
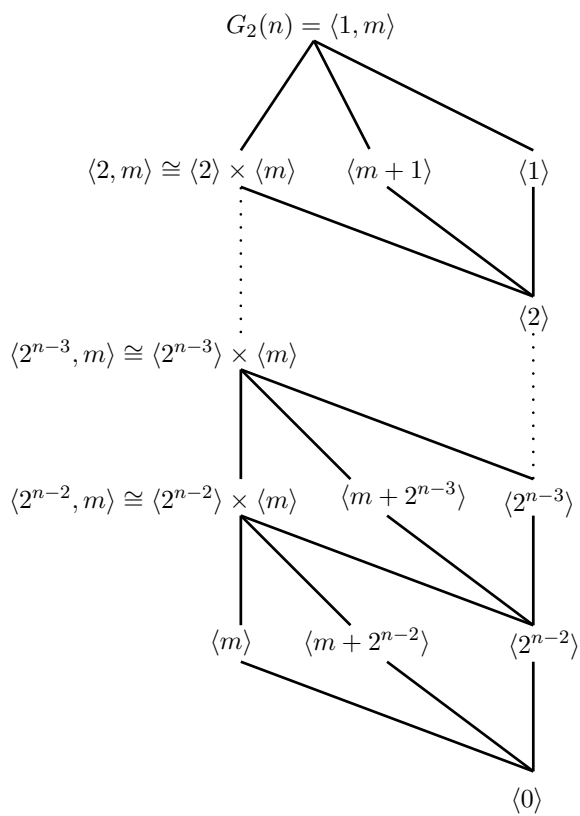
%%%%%%%%%%%%%%%%%%%%%%
\end{proof}
%%%%%%%%%%%%%%%%%%%%%%%%%%
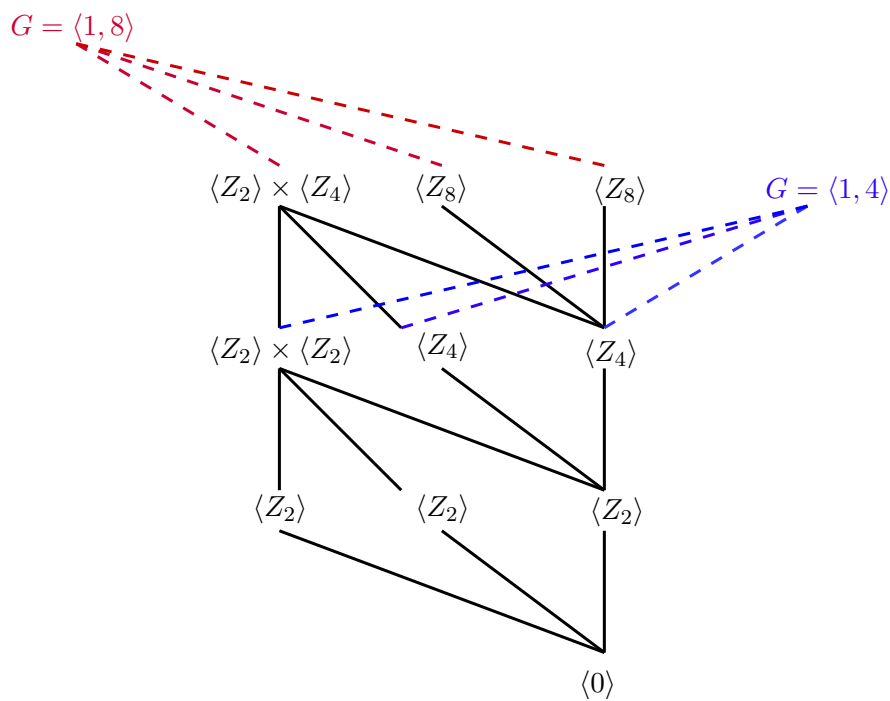
\begin{figure}[ht!]
\centering % \footnotesize
\caption{The lattice of $G_2(3)$ and $G_2(4)$ together.}\label{fig3}
~\vspace*{.5cm}\\
\scalebox{1} % Change this value to rescale the drawing.
{
\begin{pspicture}(0,-4.552969)(12.042812,4.552969)
\definecolor{color76}{rgb}{0.2,0.2,1.0}
\definecolor{color958}{rgb}{0.2,0.0,1.0}
\definecolor{color1006}{rgb}{0.8,0.0,0.0}
\definecolor{color77}{rgb}{0.8,0.0,0.2}
\psline[linewidth=0.04cm,linecolor=color77,linestyle=dashed,dash=0.16cm 0.16cm](3.5409374,2.5345314)(0.8409375,4.1545315)
\psline[linewidth=0.04cm,linecolor=color77,linestyle=dashed,dash=0.16cm 0.16cm](5.7009373,2.5345314)(0.8409375,4.1545315)
\psline[linewidth=0.04cm,linecolor=color1006,linestyle=dashed,dash=0.16cm 0.16cm](7.8609376,2.5345314)(0.8409375,4.1545315)
\psline[linewidth=0.04cm](3.5409374,1.9945313)(3.5409374,0.37453124)
\psline[linewidth=0.04cm](7.8609376,1.9945313)(7.8609376,0.37453124)
\psline[linewidth=0.04cm](7.8609376,-2.3254688)(7.8609376,-3.9454687)
\psline[linewidth=0.04cm](5.7009373,-2.3254688)(7.8609376,-3.9454687)
\psline[linewidth=0.04cm](3.5409374,-2.3254688)(7.8609376,-3.9454687)
\psline[linewidth=0.04cm](3.5409374,1.9945313)(7.8609376,0.37453124)
\psline[linewidth=0.04cm](5.7009373,1.9945313)(7.8609376,0.37453124)
\psline[linewidth=0.04cm](3.5409374,1.9945313)(5.1609373,0.37453124)

\usefont{T1}{ptm}{m}{n}
\rput(7.7723436,-4.3754687){$\langle 0\rangle$}
\usefont{T1}{ptm}{m}{n}
\rput(7.952344,0.02453125){$\langle Z_4\rangle$}
\usefont{T1}{ptm}{m}{n}
\rput(5.7123437,0.12453125){$\langle Z_4\rangle$}
\usefont{T1}{ptm}{m}{n}
\rput(8.032344,-2.0754688){$\langle Z_2\rangle$}
\usefont{T1}{ptm}{m}{n}
\rput(3.5623437,2.2245312){$\langle Z_2\rangle\times\langle Z_4\rangle$}
\usefont{T1}{ptm}{m}{n}
\rput(5.6923437,2.2245312){$\langle Z_8\rangle$}
\usefont{T1}{ptm}{m}{n}
\rput(8.072344,2.2045312){$\langle Z_8\rangle$}
\usefont{T1}{ptm}{m}{n}
\rput(0.8023437,4.364531){\color{color77}$G=\langle 1,8\rangle$}
\psline[linewidth=0.04cm](7.8609376,-0.16546875)(7.8609376,-1.7854687)
\psline[linewidth=0.04cm](3.5409374,-0.16546875)(7.8609376,-1.7854687)
\psline[linewidth=0.04cm](5.7009373,-0.16546875)(7.8609376,-1.7854687)
\psline[linewidth=0.04cm](3.5409374,-0.16546875)(5.1609373,-1.7854687)
\usefont{T1}{ptm}{m}{n}
\rput(5.7123437,-2.0354688){$\langle Z_2\rangle$}
\usefont{T1}{ptm}{m}{n}
\rput(3.5523438,-2.0354688){$\langle Z_2\rangle$}
\usefont{T1}{ptm}{m}{n}
\rput(3.5623437,0.06453125){$\langle Z_2\rangle\times\langle Z_2\rangle$}
\psline[linewidth=0.04cm](3.5409374,-0.16546875)(3.5409374,-1.7854687)
\psline[linewidth=0.04cm,linecolor=color76,linestyle=dashed,dash=0.16cm 0.16cm](10.560938,1.9945313)(7.8609376,0.37453124)
\psline[linewidth=0.04cm,linecolor=color958,linestyle=dashed,dash=0.16cm 0.16cm](10.560938,1.9945313)(5.1609373,0.37453124)
\psline[linewidth=0.04cm,linecolor=blue,linestyle=dashed,dash=0.16cm 0.16cm](3.5409374,0.37453124)(10.560938,1.9945313)
\usefont{T1}{ptm}{m}{n}
\rput(10.842343,2.2045312){\color{color958}$G=\langle 1,4\rangle$}

\end{pspicture}
}
\end{figure}

%%%%%%%%%%%%%%%%%%%%%%%%%%%%%%%%%%%%%%
%\vspace{5cm}

\vskip 3mm

%\noindent{\bf Acknowledgement.}  We are very grateful from an anonymous referee for his/her suggestions and helpful remarks.
%The research of the authors are partially supported by the University of Kashan under grant no 785149/70.
\pagebreak

\end{document}